\documentclass[12pt]{amsart} 

\usepackage{amssymb} 
\usepackage{amsmath} 
\usepackage{fullpage} 
\usepackage{url}

\theoremstyle{definition}

\theoremstyle{remark}

\def\Er{{\mathbb E}}

\def\Pr{{\mathbb P}}
\def\Qr{{\mathbb Q}}
\def\Rr{{\mathbb R}}

\def\Fc{{\mathcal{F}}}

\def\Sc{{\mathcal{S}}}

\def\one{{\rm \bf 1}}

\def\essinf{\operatorname{ess.\!inf}}

\def\({\left(}     
\def\){\right)}    
\def\[{\left[}     
\def\]{\right]}

\begin{document}
\title{Remark on the Paper ``Entropic Value-at-Risk: A New Coherent Risk Measure" by Amir Ahmadi-Javid, J. Opt. Theory and Appl., 155 (2001),1105--1123  \footnote{{\it AMS-}Classification 90B50, 91B06, 91B16, 91G99} } 
\author{Freddy Delbaen}
\address{Departement f\"ur Mathematik, ETH Z\"urich, R\"{a}mistrasse
   101, 8092 Z\"{u}rich, Switzerland}
 \address{Institut f\"ur Mathematik,
 Universit\"at Z\"urich, Winterthurerstrasse 190,
 8057 Z\"urich, Switzerland}
\date{First version November 6, 2013, this version \today}

\begin{abstract}
The paper mentioned in the title introduces the entropic value at risk.  I give some extra comments and using the general theory make a relation with some commonotone risk measures.\end{abstract}

\maketitle

\section{Introduction}

In \cite{AAJ}, see also \cite{AAJ_ad} for corrections and precisions, Amir Ahmadi-Javid introduces the Entropic Value at Risk (EVAR). He also relates this risk measure to the CVaR or expected shortfall measure and gives some generalisations.  In this note I make some remarks, relating the EVAR to some other commonotone risk measure.  In particular I show that there is a best commonotone risk measure that is dominated by EVAR.  The presentation uses the concept of monetary utility functions, which has the advantage that it has a better economic interpretation.  Up to some sign changes it is the same as in the paper of \cite{AAJ}.  The general theory of such coherent utility functions was introduced in \cite{ADEH1}, \cite{ADEH2}, \cite{Del}. For a more complete presentation I refer to \cite{Pisa} and \cite{FDbook}.

We make the standing assumption that all random variables are defined on an atomless probability space $(\Omega,\Fc,\Pr)$.  For a level $\alpha, 0<\alpha<1$, the Entropic Value at Risk of a bounded random variable, $\xi$, is defined as
$$
e_\alpha(\xi)=\sup_{0<z}\(-\frac{1}{z}\log \Er\[\frac{\exp(-z\xi)}{\alpha}\]\).
$$
The reader will --- already said above --- notice that this is up to sign changes, the definition of \cite{AAJ}.  The entropic value at risk is the supremum of entropic utility functions perturbed by the introduction of the level $\alpha$.  For $\alpha=1$ we would find the $\Er_\Pr[\xi]$. For $\alpha\rightarrow 0$, we get $\essinf \xi$. The mapping
$$
e_\alpha\colon L^\infty\rightarrow \Rr,
$$
defines a coherent utility function, see \cite{AAJ}.  Using duality arguments,  it is shown there that
$$
e_\alpha(\xi)=\inf\left\{\Er_\Qr[\xi]\mid H(\Qr|\Pr)=\Er\[\frac{d\Qr}{d\Pr}\log\(\frac{d\Qr}{d\Pr}\)\]\le -\log(\alpha)\right\}.
$$
The proof of the equality uses a standard duality argument (see the last line on page 1111 of the paper) but also needs the (not mentioned) property that the sets $\{\Qr\mid H(\Qr\mid\Pr)\le \beta\}$ ($\beta\ge 0$) are weakly compact in $L^1$. For later use let us denote
$$
\Sc= \{\Qr\mid H(\Qr\mid\Pr)\le -\log(\alpha)\}
$$
In the next sections we will calculate the values for indicators and using these we will prove that $e_\alpha$ is not commonotone.  We will show that $e_\alpha$ is between two commonotone utility functions.
\section{The value for indicators}

In case $\xi=\one_A$ is the indicator of a set $A\in\Fc$, we can more or less calculate the value $e_\alpha(\one_A)$. Because the function $x\log(x)$ is convex, the entropy decreases if we replace $\Qr$ by the probability given by the density $\Er\[\frac{d\Qr}{d\Pr}\mid A,A^c\]$.  Because the space is atomless, it is easily seen that the value
$$ 
e_\alpha(\one_A)=\inf\left\{\Er_\Qr[\one_A]\mid H(\Qr|\Pr)=\Er\[\frac{d\Qr}{d\Pr}\log\(\frac{d\Qr}{d\Pr}\)\]\le -\log(\alpha)\right\},
$$
is given by a minimum attained by a function of the form ($a=\Pr[A]$):
$$
\frac{d\Qr}{d\Pr}=\lambda \frac{\one_A}{a} +(1-\lambda)\frac{\one_{A^c}}{1-a}.
$$
In this case the value $\Er_\Qr[\one_A]$ is simply $\lambda$. The entropy of such a measure is given by
$$F(\lambda,a)=\lambda\log(\lambda)+(1-\lambda)\log(1-\lambda) - \lambda\log(a)-(1-\lambda)\log(1-a).$$
The  function $F$  takes a unique minimum equal to zero at $\lambda=a$ and for $\lambda=0$ it gives the value $-\log(1-a)$, whereas for $\lambda=1$ we get $-\log(a)$.  The derivatives of $F$ satisfy:
$$\frac{\partial F}{ \partial \lambda}=\log\(\frac{\lambda}{1-\lambda}\) -\log\(\frac{a}{1-a}\)<0\quad\quad  \frac{\partial F}{ \partial a}=\frac{a-\lambda}{a(1-a)}>0$$
$$
\frac{\partial^2 F}{ \partial \lambda^2}=\frac{1}{\lambda}+\frac{1}{1-\lambda};\,\,
\frac{\partial^2 F}{ \partial \lambda\partial a}=-\frac{1}{a}-\frac{1}{1-a};\,\,
\frac{\partial^2 F}{ \partial a^2}=\frac{\lambda}{a^2}+\frac{1-\lambda}{(1-a)^2}
$$
for $0<\lambda < a$, $1-\alpha<a<1$.  It follows easily that for $a\le (1-\alpha)$, the minimum in the expression for $e_\alpha(\one_A)$ is attained for $\lambda=0$ and consequently $e_\alpha(\one_A)=0$.  For $1> a >1-\alpha$ there is one solution of the equation 
$$F(\lambda,a)=\lambda\log(\lambda)+(1-\lambda)\log(1-\lambda) - \lambda\log(a)-(1-\lambda)\log(1-a)=-\log(\alpha)$$ 
that is smaller than $a$. The function $\Lambda(a)=e_\alpha(\one_A) < a$ is well defined as an implicit function.  From the definition of $e_\alpha$, it is already seen that $\Lambda(a)$ must be increasing but the implicit function theorem gives that $\frac{d\Lambda}{da}>0$ on $1-\alpha<a<1$. For $a\rightarrow 1$, the derivative tends to $+\infty$.  The Hessian of $F$ is positive definite and hence on $1-\alpha<a<1,0<\lambda < a$, the function $F$ is strictly convex. If we take the derivative of $\frac{\partial F}{ \partial a} + \frac{\partial F}{ \partial \lambda}\frac{d\Lambda}{da}=0$ with respect to $a$ we find at the points $(a,\Lambda(a))$ that
$$
\frac{\partial^2 F}{ \partial a^2}+2\,\frac{\partial^2 F}{ \partial \lambda\partial a}\frac{d\Lambda}{da} +
\frac{\partial^2 F}{ \partial \lambda^2}\(\frac{d\Lambda}{da}\)^2 = - \frac{\partial F}{ \partial \lambda}\frac{d^2\Lambda}{da^2}.
$$
Since the Hessian of $F$ is positive definite the left hand side is always positive. Hence we must have $\frac{d^2\Lambda}{da^2}>0$ on $1-\alpha<a<1,0<\lambda < a$. The function $\Lambda$ is therefore strictly convex on $1-\alpha\le a \le 1$.  Summarising $\Lambda\colon[0,1]\rightarrow[0,1]$ is convex, $\Lambda(0)=0,\Lambda(1)=1$ and $\Lambda(a)=0$ for $0\le a\le 1-\alpha$.
\section{Relation with commonotone utilities}
For a convex function $f\colon [0,1]\rightarrow [0,1]$ with $f(0)=0,f(1)=1$ we can define a commonotone utility.  The relation with convex games and non-additive expectations is well known, see \cite{Schm}, \cite{FDbook}, \cite{Pisa}. The basic ingredient is the scenario set defined as
$$
\Sc_f=\{\Qr\mid \text{ for all } A\in \Fc: \Qr[A]\ge f(\Pr[A])\}.
$$
With this we associate the utility
$$
u_f(\xi)=\inf\{\Er_\Qr[\xi]\mid  \Qr\in\Sc_f\}.
$$
We immediately get $\Sc\subset \Sc_\Lambda$ and $\Sc_c\subset \Sc$, where $c$ denotes the convex function defined as $c(x)=0$ for $x\le 1-\alpha$ and $c(x)=\frac{x-1+\alpha}{\alpha}$ for $1-\alpha\le x\le 1$.  The latter defines the utility function known as CVAR at level $\alpha$, see \cite{Pisa} or \cite{FDbook}.  We get that
$$
CVAR_\alpha \ge e_\alpha \ge u_\Lambda.
$$
The calculation of $u_\Lambda$ can be done using Ryff's theorem (see \cite{Ryff}, \cite{FDbook}). We get:
$$
u_\Lambda(\xi)=\int_0^1 q_{1-a} \frac{d\Lambda}{da} da = \int_{1-\alpha}^1 q_{1-a} \frac{d\Lambda}{da} da,
$$
where $q$ is a quantile function of $\xi$.

For indicators $\one_A$ we have $u_\Lambda(\one_A)=e_\alpha(\one_A)$. Therefore we get that $u_\Lambda$ is the greatest commonotone utility that is dominated by $e_\alpha$.  The utility $e_\alpha$ is not commonotone.  Indeed if it were, then for $1-\alpha<a < b < 1$, $A\subset B$, $a=\Pr[A], b=\Pr[B]$ we would find an element of $\Sc$, say $\Qr$ such that $\Er_\Qr[\one_A]= \Lambda (a), \Er_\Qr[\one_B]=\Lambda(b)$.  However the minimisers for indicators have a special form and this would imply that $\frac{\Lambda(a)}{a}=\frac{\Lambda(b)}{b}$ as well as $\frac{\Lambda(b)}{b}=\frac{1-\Lambda(a)}{1-a}$.  The only solution would then be $\Lambda(a) = a$, a contradiction.  The reader can check that for $\xi=\one_A+\one_B$ we have that $e_\alpha(\xi) > u_\Lambda(\xi)$.  With some extra effort the reader can also see that for random variables $\xi$ such that the distribution has a support of more than two points, we have the same strict inequality.
\section{The Kusuoka Representation}
 In \cite{kusuoka}, Kusuoka proved that coherent utilities that only depend on the law of the random variables are given by a convex set $K$  of probabilities on $[0,1]$. They are in fact averages of CVAR utilities: 
 $$
 u(\xi) = \inf \left\{ \int_{[0,1]} CVAR_x(\xi)\, \nu(dx)\mid \nu \in K\right\}.
 $$
To find the set $K$ needed to represent $e_\alpha$, we proceed as in \cite{FDbook}.  Unfortunately besides a transformation of the problem we cannot give a more explicit analytic form for the set $K$.  The first step is to find decreasing probability densities $\eta$ on $[0,1]$ such that $\int_0^1 \eta(x)\log(\eta(x))\,dx \le -\log (\alpha)$.  Then we write each $\eta$ in the form $\eta(x)=\int_{(x,1]}\frac{1}{a}\,\nu(da)$ where $\nu$ is a probability on $[0,1]$.  The set $K$ is then the set of the measures $\nu$ so obtained. As shown in \cite{FDbook} the set $K$ is a weak$^*$ compact convex set of probability  measures on $(0,1]$.  For  more information the reader could also check Remark 41 on page 95 and exercise 18 in the same book.  Somewhat more elegant descriptions can  be given but essentially no better characterisation was found. I do not pursue this direction here.
\section{Acknowledgement}
The author thanks dr.\! \!Marcus Wunsch, UBS, Switzerland for drawing his attention to this paper and for valuable comments.  The present paper was written while the author was on visit at Fudan University in Shanghai and at Shandong University in Jinan.  The author thanks both institutions for their hospitality and for discussions on the topic.

\end{document}